\documentclass{amsart}
\usepackage{amsfonts}
\usepackage{amsmath}
\usepackage{amssymb}
\usepackage{graphicx,array}

\newtheorem{theorem}{Theorem}[section]
\theoremstyle{plain}

\newtheorem{corollary}[theorem]{Corollary}
\newtheorem{lemma}[theorem]{Lemma}

\newcommand\straightparen[1]{{\rm(}{#1}{\rm\,)}}
\newcommand\usfrac[2]{{#1}/{#2}}

\def\dsum{\displaystyle\sum}
\newcommand{\nnn}[1]{{\rm(}\ref{#1}{\rm)}}

\numberwithin{equation}{section}

\numberwithin{equation}{section}

\newcommand{\fa}{\mathfrak{a}}

\newcommand{\gD}{\Delta}

\newcommand{\R}{\mathbb{R}}
\newcommand{\C}{\mathbb{C}}
\newcommand{\N}{\mathbb{N}}
\newcommand{\Z}{\mathbb{Z}}
\newcommand{\bc}{\mathbf{c}}

\newcommand{\Supp}{\mathrm{Supp}}

\newcommand{\e}{\varepsilon}
\renewcommand{\l}{\lambda}

\newcommand{\inner}[2]{\langle #1,#2\rangle}

\begin{document}
\makeatletter
\title[Paley-Wiener Theorem]{The Paley-Wiener Theorem for the Jacobi transform\\ and the local Huygens'
principle\\ for root systems with even multiplicities}
\author{Thomas Branson}
\address{Department of Mathematics, University of Iowa,
Iowa City IA 52242 USA}
\email{thomas-branson@uiowa.edu}
\urladdr{http://www.math.uiowa.edu/\symbol{126}branson}
\author{Gestur \'{O}lafsson}
\address{Department of Mathematics, Louisiana State University,
Baton Rouge LA 70803, USA}
\email{olafsson@math.lsu.edu}
\urladdr{http://www.math.lsu.edu/\symbol{126}olafsson}
\author{Angela Pasquale}
\address{Laboratoire et D\'epartement de Math\'ematiques,
  Universit\'e de Metz, France}
\email{pasquale@math.univ-metz.fr}
\urladdr{http://www.math.univ-metz.fr/\symbol{126}pasquale}
\thanks{GO was supported by NSF grants  DMS-0139783,
and DMS-0402068, and
by the DFG-Schwerpunkt \textit{Global Methods in Complex Geometry}.
TB was partially supported by the Erwin Schr\"odinger Institute.}
\date{August 3, 2005}
\dedicatory{Dedicated to Gerrit van Dijk on the occasion of his 65th birthday}

\begin{abstract}
This note is a continuation of the previous paper \cite{BOP1} by same authors.
Its purpose is to extend the results of \cite{BOP1} to the context of
root systems with even multiplicities. Under the even multiplicity
assumption, we prove  a local Paley-Wiener theorem for the
Jacobi transform and the strong Huygens' principle for the wave
equation associated with the modified compact Laplace operator.
\end{abstract}
\maketitle
\makeatother

\section*{Introduction}

\noindent
Harish-Chandra's theory of spherical functions on Riemannian symmetric
spaces of the noncompact type, resp.\ of the compact type, was
generalized in the late 1980s by G.\ Heckman and E.\ Opdam into
the theory of hypergeometric functions, resp.\ of Jacobi polynomials,
associated with root systems. Correspondingly, the noncompact and compact
spherical transforms have as natural generalizations
the Opdam's hypergeometric transform and the Jacobi transform.
We refer to \cite{HS94}, \cite{OpdActa}, \cite{Opd00} and references therein for more information. It is then a natural question, how known results in harmonic analysis on symmetric spaces can be extended
to these new integral transforms.

In \cite{BOP1} we proved a Paley-Wiener theorem for
the spherical transform on compact symmetric spaces with
even multiplicities. The Paley-Wiener theorem was then
used to show (by three different methods) that the local
strong Huygens' principle holds for the modified wave equation on these spaces.
In this note we show how our results can be extended
to prove a local Paley-Wiener theorem for the Jacobi transform
and then the strong Huygens principle for the modified
Laplacian associated with a root system with even multiplicities.

In the context of special functions associated with root systems,
the symmetric space is replaced by a triple $(\mathfrak a, \Delta, m)$,
consisting of an $n$-dimensional real Euclidean vector space $\fa$, a root
system $\Delta$ in the dual $\fa^*$ of $\fa$, and a multiplicity
function $m$ on $\Sigma$. Hence $\alpha \mapsto m_\alpha$ is a map
on $\Delta$ which is invariant with respect to the Weyl group $W$ of
$\Delta$. We will write $m\geq 0$ to indicate that
$m_\alpha\geq 0$ for all $\alpha \in \Delta$. The two results
in the title will be proven under the additional assumption that
$\Delta$ is reduced, i.e. $2\alpha \notin \Delta$ for all $\alpha \in
\Delta$, and that $m$ is even, i.e. $m_\alpha \in 2\N$ for all
$\alpha \in \Delta$.
In the geometric situation of symmetric spaces, even multiplicity functions
occur only on reduced root systems. In our more general context, the two
conditions have to be imposed.
We remark that our notation
is based on the theory of symmetric spaces which differs from the  Heckman-Opdam
notation in the following ways. The root system $R$ used by Heckman and Opdam
is related to our root system $\Delta$ by $R=\{2\alpha\mid \alpha \in \Delta\}$
and the multiplicity function $k$ in Heckman-Opdam's work
is given by $k_{2\alpha}=m_\alpha/2$.

Compared to \cite{BOP1}, this note contains two new ingredients: an
explicit formula for the Jacobi polynomials in the even multiplicity case,
and a detailed study of the convergence of certain Jacobi expansions. It
is also important to notice that in \cite{BOP1} several results have
were proven under the more restrictive condition (always satisfied by Riemannian
symmetric spaces with even multiplicities) that
the multiplicity function is a constant function, i.e.,
there is $m\in \N$ such that $m_\alpha=2m$ for all $\alpha \in \Delta$.
This condition will be dropped in the present note.
We also remark that we will not work out all details. Our
aim is to prove what is necessary to apply the results
and methods from \cite{BOP1}.

\section{Jacobi polynomials and the Jacobi Transform}
\noindent
Let $P$ be the lattice of restricted weights. This consists of
all elements $\mu \in \fa_\C^*$ so that
$\mu_\alpha:=\usfrac{\inner{\mu}{\alpha}}{\inner{\alpha}{\alpha}} \in \Z$
for all $\alpha \in \Delta$.
In particular $P \subset \fa^*$. Then $A_\C:={\rm Hom}_\Z(P,\C^\times)$
is a complex torus with Lie algebra $\fa_\C$. It admits
the decomposition $A_\C=AT$ with $A:=\exp(\fa)$ and $T=\exp(i\fa)$.
Recall that harmonic analysis of $K$-invariant function spaces
on $G/K$, respectively $U/K$, can be reduced to harmonic analysis
of Weyl group invariant objects on a maximal flat subspace. Therefore
the Lie groups $A$ and $T$ are respectively replacements for the
symmetric spaces $G/K$ and $U/K$.

\subsection{Jacobi polynomials}
Let $\C[A_\C]$ denote the space of finite $\C$-linear
combinations of elements $e^{\mu}$ with $\mu \in P$, and
let $\C[A_\C]^W$ be the subspace of $W$-invariant elements. We
will always assume that $m\geq 0$.
Fix a set of positive roots $\Delta^+$ in $\Delta$ and
set $\delta(m,t):=\prod_{\alpha \in \Delta^+}
|e^{\alpha(\log t)}-e^{-\alpha(\log t)}|^{m_\alpha}$.
Define an inner product $\inner{\cdot}{\cdot}_m$ on $\C[A_\C]^W$ by
$$\inner{f}{g}_m:=\int_T f(t) \overline{g(t)} \; \delta(m,t) \,dt\,,$$
where $dt$ is the normalized Haar measure on $T$.

Let $P^+:=\{\mu \in \fa^* \mid \text{$\mu_\alpha\in \Z^+$ for
all $\alpha \in \Delta^+$}\}$.
For $\mu \in P^+$ and $t \in T$ let $P(m,\mu,t)$ denote
the Jacobi polynomial with spectral parameter $\mu$. This is
defined as follows. The orbit sums
$$
M(\mu):=\sum_{\nu \in W\mu} e^\nu
$$
form a basis of $\C[A_\C]^W$ as $\mu$ varies in $P^+$ because each
$W$-orbit in $P$ intersects $P^+$ in exactly one point. The
Jacobi polynomial $P(m,\mu)$ is the exponential polynomial
\begin{equation}\label{eq:Jacobipoly}
P(m,\mu):=\sum_{\mu\ge\nu \in P^+} c_{\mu\nu}(m) M(\nu)
\end{equation}
where the coefficients $c_{\mu\nu}(m)$ are defined by the following
conditions, cf.\ \cite{HS94}, \S 1.3:
\begin{enumerate}
\renewcommand{\theenumi}{\roman{enumi}}
\renewcommand{\labelenumi}{(\theenumi)}
\item
$c_{\mu\mu}(m)=1$\,;
\item
$\inner{P(m,\mu)}{M(\nu)}_m=0$ for all $\nu \in P^+$ with
$\nu < \mu$.
\end{enumerate}
Observe that by definition,  $P(m,\mu,t)$ extend holomorphically
to $A_\C$ as a function of $t$.
Moreover $\{P(m,\mu) \mid \mu \in P^+\}$ is a basis for $\C[A_\C]^W$ which is orthogonal
with respect to the inner product
$\inner{\cdot}{\cdot}_m$ (cf.  \cite{HS94}, Corollary 1.3.13).

The $L^2$-norm of $P(m,\mu)$ is expressed in terms of two $\bc$-functions. The first
one, denoted $\bc(m,\lambda)$, is the $\mathbf{c}$-function
defined by means of the Gindikin-Karpelevic formula as in \cite{BOP1}, (2.10).
In particular, for reduced root systems and even multiplicity functions,
\begin{equation}\label{eq:c}
\frac{1}{\mathbf{c}(m,\lambda)}=C\prod_{\alpha\in\Delta^{+}}
\prod_{k=0}^{m_\alpha/2-1}(\lambda_{\alpha}+k)
\end{equation}
where the constant $C$ is given by
\begin{equation} \label{eq:C}
C=\prod_{\alpha\in\Delta^+}\prod_{k=0}^{m_\alpha/2-1}\frac{1}{\rho_\alpha+k}
\end{equation}
and
$\rho:=\rho(m):=\frac{1}{2}\sum_{\alpha\in\Delta^+} m_\alpha \alpha$
as usual.

The second function, denoted $\bc^*(m,\l)$, is a dual $\bc$-function.
In the case of reduced root systems and even multiplicities, these two functions are
related by the formula
$$\bc^*(m,\l)=C\, (-1)^{\sum_{\alpha\in\Delta^+} m_\alpha/2}
  \bc(m,\lambda)\,$$
with $C$ as in (\ref{eq:C}).
See \cite{HS94} (3.5.2) for the general definition of $\bc^*$.
Again assuming even multiplicities, one has for $\mu \in P^+$
(see \cite{HS94}, Corollary 3.5.3):
\begin{equation}
\label{eq:normP}
\|P(m,\mu)\|^2_m=\dfrac{|W|}{C}\,
\dfrac{\bc^*(m,-\mu-\rho)}{\bc(m,\mu+\rho)}
=  |W| \prod_{\alpha\in\Delta^+}\prod_{k=0}^{m_\alpha/2-1}
\dfrac{(\mu_{\alpha}+\rho_\alpha+k)}
{(\mu_{\alpha}+\rho_\alpha-k)}\,.
\end{equation}
The positivity of the last member in (\ref{eq:normP}) depends on
the following lemma.

\begin{lemma}\label{le-dim}
Assume that $\Delta$ is a reduced root system and that for all $\alpha\in\gD$
the multiplicities $m_\alpha$
are even.
Then $\rho_\alpha\geq m_\alpha/2$ for every $\alpha \in\Delta^+$.
\end{lemma}
\begin{proof}
Observe first that
\begin{equation}\label{eq:rhoalphapositive}
\rho_\alpha>0 \quad \text{for every $\alpha \in\Delta^+$}.
\end{equation}
Indeed, if $\alpha$ is a simple root in a reduced root system $\Delta$, then $\rho_\alpha=m_\alpha/2$. See \cite{BOP1}, Lemma 2.7, (1).

We prove the inequalities in the statement by induction on the length $\ell(\alpha)$ of a root $\alpha \in\Delta^+$. Recall that, if $\Pi:=\{\alpha_1,\dots,\alpha_n\}$ is
the basis of simple roots in
$\Delta^+$ and $\alpha=\sum_{j=1}^n n_j \alpha_j$, then $\ell(\alpha):=\sum_{j=1}^n n_j$.
The inequalities being true when $\ell(\alpha)=1$, we then  suppose that $\ell(\alpha)>1$.
Because of (\ref{eq:rhoalphapositive}), there must be $\gamma \in \Pi$ so that
$\inner{\gamma}{\alpha}> 0$.
Since $\ell(\alpha)>1$, the root $\alpha$ is not simple. Hence
$s_\gamma \alpha:=\alpha-2\frac{\inner{\gamma}{\alpha}}{\inner{\gamma}{\gamma}}\; \gamma$
satisfies $s_\gamma \alpha\in \Delta^+$ and $\ell(s_\gamma \alpha)<\ell(\alpha)$.
We can therefore apply the inductive hypothesis to $s_\gamma \alpha$.
Noticing that $s_\gamma \rho=\rho-m_\gamma \gamma$, we get:
$$
\begin{array}{l}
\dfrac{\inner{\rho}{\alpha}}{\inner{\alpha}{\alpha}}=
\dfrac{\inner{s_\gamma\rho}{s_\gamma\alpha}}{\inner{s_\gamma\alpha}{s_\gamma\alpha}}
=\dfrac{\inner{\rho}{s_\gamma\alpha}}{\inner{s_\gamma\alpha}{s_\gamma\alpha}}-
m_\gamma\; \dfrac{\inner{\gamma}{s_\gamma\alpha}}{\inner{s_\gamma\alpha}{s_\gamma\alpha}}\\
\qquad\geq m_{s_\gamma \alpha}/2+ m_\gamma\; \dfrac{\inner{\gamma}{\alpha}}{\inner{s_\gamma\alpha}{s_\gamma\alpha}}
\geq m_{s_\gamma\alpha}/2=m_{\alpha}/2.
\end{array}
$$
\end{proof}

\subsection{Hypergeometric functions}
Let $F(m,\l,a)$ denote the hypergeometric function of spectral
parameter $\l \in \fa_\C^*$. By definition,
this is the unique solution of the Heckman-Opdam
hypergeometric system of differential equations corresponding
to the spectral parameter $\l$ which is analytic in the space
parameter $a\in A$ and normalized by the condition $F(m,\l,e)=1$.
Here $e$ denotes the identity element of $A$. See e.g.\ \cite{HS94}, \S 4.
If the triple $(\mathfrak a, \Delta, m)$ is geometric (i.e.\ arises from a
Riemannian symmetric space of the noncompact type $G/K$), then $F(m,\l,a)$
agrees with the restriction to $A=\exp(\mathfrak a)$ of the
spherical function $\varphi_\l$ from \cite{BOP1}, (2.7).

\subsection{Explicit formulas}
The generalization of Lemma 2.5 in  \cite{BOP1} and
the Weyl group invariance stated in  Theorem 2.11 in \cite{BOP1}
is given as follows.

\begin{lemma} Suppose that $m \geq 0$ is a fixed multiplicity function.
Then
\begin{equation} \label{eq:restrictiontoA}
F(m,\mu+\rho,a)=\mathbf{c}(m,\mu+\rho)P(m,\mu,a)\,.
\end{equation}
Furthermore,
\begin{equation}
c(m, w(\mu-\rho)-\rho) P(m, w(\mu-\rho)-\rho)=c(m,\mu)P(m,\mu)
\end{equation}
for all $\mu \in P^+$ and $w\in W$.
\end{lemma}
\begin{proof}
See \cite{HS94}, (4.4.10).
\end{proof}
Formula (\ref{eq:restrictiontoA}) provides a holomorphic extension
of the hypergeometric function $F(m,\mu+\rho)$ to $A_\C$. Moreover,
the second relation shows that we can extend
the definition of $P(m,\mu)$ to $\mu \in P$.
Observe, however, that in the context of special functions associated
with root systems, the functions $F(m,\l)$ and $P(m,\mu)$ do not
generally admit integral representations. Likewise, these functions
cannot be considered as matrix coefficients of group representations. Nevertheless,
the representation dimensions, or Plancherel density, $d(\mu)=d(m,\mu)$ can still be introduced
by means of Vretare's formula:
\begin{equation}\label{eq:d}
d(m,\mu):=\lim_{\varepsilon\to 0}\;
\dfrac{\mathbf{c}(m,-\rho+\varepsilon)}{\mathbf{c}(m,\mu+ \rho)\mathbf{c}(m,-\mu+\rho+\varepsilon)}  \,.
\end{equation}
As in \cite{BOP1}, this formula simplifies because of Lemma \ref{le-dim},
which allows us to compute the limit in (\ref{eq:d}) as the quotient
of the limits of the $\bc$-functions appearing in the numerator and in the denominator.

\begin{corollary}\label{cor:d}
Assume that $\Delta$ is reduced and all $m_\alpha$ are even for all $\alpha\in\gD$.
Then the following properties hold:
\begin{enumerate}
\item For all $\mu\in P^+$ we have
$$
d(m,\mu)=\frac{\bc(m,-\rho)}{\bc(m,\mu+\rho)\bc(m,-(\mu+\rho))}\,.
$$
\item The function $d(m,\mu)$ extends as to a polynomial function on $\fa^*_\C$ given by
\begin{equation*}
d(m,\lambda )
=\prod_{\alpha\in \Delta^{+}}\prod_{k=0}^{m_\alpha/2-1}
\frac{k^{2}-(\lambda +\rho)_{\alpha}^{2}}{k^2-\rho_\alpha^2}\, .
\end{equation*}
\end{enumerate}
\end{corollary}

Assuming that $\Delta$ is reduced and all multiplicities are even, the
function $\delta(m,t)$ extends
to a $W$-invariant holomorphic function on $A_\C$. Moreover, because of
Theorem 5.1(c) in \cite{OP04}, there exists a
$W$-invariant differential operator $D$ on $A$ with coefficients
which are holomorphic on $A_\C$ so that for all $\lambda \in
\mathfrak a_\C^*$ and all $a\in A$ we have
\begin{equation}\label{eq:formulaF}
\delta(m,a)F(m,\lambda,a)=\dfrac{1}{d(m,\lambda-\rho)}
\, D\left(\sum_{w\in W} a^{w\lambda}\right).
\end{equation}
The right-hand side of (\ref{eq:formulaF}) is holomorphic in
$a\in A_\C$ and therefore a provides holomorphic extension of
$\delta(m,a)F(m,\lambda,a)$. It also follows from (\ref{eq:formulaF}) that
\begin{equation}
\label{eq:Dzero}
D\left(\sum_{w\in W} a^{w\lambda}\right)=0
\end{equation}
for all $\lambda\in \mathfrak a_\C^*$ satisfying $\lambda_\alpha\in
\pm \{0,1,\dots, m_\alpha/2-1\}$ for some $\alpha \in \Delta^+$.

As an easy corollary of the above, we obtain the following
formulas for the Jacobi polynomials.

\begin{corollary}\label{cor:formulaJacobi}
Suppose that the root system is reduced and all multiplicities are even.
Then for all $\mu \in P$
and all $a\in A_\C$ we have
$$\delta(m,a)P(m,\mu,a)=\dfrac{1}{\bc(m,\mu+\rho)d(m,\mu)}
\, D\left(\dsum_{w\in W} a^{w(\mu+\rho)}\right)\,$$
and
\begin{equation} \label{eq:formulanormJacobi}
\begin{array}{rl}\delta(m,a)\dfrac{P(m,\mu,a)}{\|P(m,\mu)\|_m^2}&=
\dfrac{|W|}{C} \; \dfrac{1}{\bc^*(m,-\mu-\rho)d(m,\mu)}
\, D\left(\dsum_{w\in W} a^{w(\mu+\rho)}\right) \\
&=\widetilde{C} \; \bc(m,\mu+\rho)
D\left(\dsum_{w\in W} a^{w(\mu+\rho)}\right),
\end{array}
\end{equation}
where
\begin{equation}
\widetilde{C}=\label{eq:Ctilde}
\dfrac{|W|}{C^2} \; \prod_{\alpha\in \Delta^+} \prod_{k=0}^{m_\alpha/2-1}(\rho_\alpha-k)\,.
\end{equation}
\end{corollary}

\subsection{Estimates for the Jacobi polynomials}
Estimates for the Jacobi polynomials and their derivatives can be deduced from (\ref{eq:restrictiontoA}), from the estimate
\begin{equation}\label{eq:estc}
\dfrac{1}{|\bc(m,\lambda)|}\leq {\rm const.} (1+\|\lambda\|)^{\sum_{\alpha\in\Delta^+} m_\alpha/2}
\end{equation}
for the $\bc$-function, and from Opdam's estimates for the hypergeometric function and its derivatives. We note that the latter estimates are only valid inside
$\{X+i H\in\mathfrak{a}_\C\mid \text{$|\alpha (H)|\le \pi /2$ for all
$\alpha\in\Delta$}\}$,
which is smaller than the domain used in \cite{BOP1}.

\begin{lemma}\label{lemma:estJacobi}
Suppose that $m\geq 0$ is a fixed multiplicity function.
\begin{enumerate}
\item
Let $t=\exp(iH)\in T$ satisfy $|\alpha(H)|\leq \pi/2$ for all $\alpha \in \Delta$. Then we have for all $\mu \in P^+$
\begin{align}\label{eq:estP}
|P(m,\mu,t)|&\leq C_1 |W|^{1/2}  (1+\|\mu\|)^{\sum_{\alpha\in\Delta^+} m_\alpha/2}
e^{\max_{w\in W} w\rho(H)} \notag\\
&\leq C_2 (1+\|\mu\|)^{\sum_{\alpha\in\Delta^+} m_\alpha/2}\,,
\end{align}
where $C_1$ and $C_2$ are some positive constants.
\item Let $I=(i_1,\dots,i_n)$ be a multi-index and set $|I|=\sum_{j=1}^n
i_j$. Let
$\partial^I=\partial(H_1)^{i_1} \cdots  \partial(H_1)^{i_n}$ be the corresponding partial differential operator associated with an orthonormal basis $\{H_1,\dots, H_n\}$ of $\mathfrak t$. Let $t=\exp(iH)\in T$ satisfy $|\alpha(H)|\leq \pi/2$ for all $\alpha \in \Delta$. Then there is a constant $C_3(t)$, depending on $t$, such that for all $\mu \in P^+$ we have
\begin{equation*}
|\partial^I P(m,\mu,t)|\leq C_3(t) (1+\|\mu\|)^{|I|+ \sum_{\alpha\in\Delta^+} m_\alpha/2}.
\end{equation*}
\end{enumerate}
\end{lemma}
\begin{proof}
The first estimate is an immediate consequence of \cite{OpdActa}, Theorem 3.15 and Proposition 6.1, together with (\ref{eq:estc}).

To prove the estimates for the derivatives, we proceed as \cite{OpdActa}, Corollary 6.2. Let $t=\exp(iH)$ be fixed as in the statement. Choose $\delta=\delta(t)>0$ so that
$T_\mu:=\{z\in \mathfrak t_\C\mid |z_j-iH_j|<\delta/\|\mu\|\} \subset
\{z=iX\in \mathfrak t\mid \text{$|\alpha(X)|<\pi/2$ for all $\alpha \in \Delta^+$}\}$ for all $\mu \in
P^+$ with $\mu \neq 0$. Cauchy's integral formula then gives
$$\partial^I P(m,\mu,t)= \int_{iH+T_\mu}
\dfrac{P(m,\lambda,z)}{(z-iH)^{I+(1,\dots,1)}}\, dz\,.$$
The result follows then easily from (\ref{eq:estP}).
\end{proof}

\subsection{The Jacobi transform}
The Jacobi transform of $f \in L^2(T)^W$ is the function $\widehat{f}(m,\cdot):P^+ \to \C$ defined by
$$\widehat{f}(m,\mu):=\inner{f}{P(m,\mu)}_m
=\int_T f(t)P(m,\mu,t^{-1})\, \delta(m,t)\; dt\,.$$
Here we have used the property that $\overline{P(m,\mu,t)}=P(m,\mu,t^{-1})$.
The inversion formula is given by
\begin{equation}
\label{eq:invJacobi}
f=\sum_{\mu \in P^+} \widehat{f}(m,\mu) \dfrac{P(m,\mu)}{\|P(m,\mu)\|^2_m}
\end{equation}
with convergence in the sense of $L^2$.
We remark that $\delta(m,t)=\delta(m,t^{-1})$ in the even multiplicity case. If $f \in C^\infty(T)^W$ we therefore obtain from Corollary
\ref{cor:formulaJacobi} that
\begin{equation}
\label{eq:Jacobiexpl}
\bc(m,\mu+\rho)d(m,\mu) \widehat{f}(m,\mu)=|W| \,
\int_T (D^*f)(t) t^{-(\mu+\rho)}\; dt,
\end{equation}
where $D^*$ denotes the formal adjoint of the differential operator $D$ with respect to the measure $dt$.

\section{The Paley-Wiener theorem for the Jacobi transform}
\noindent
In the following we shall assume that a constant $R>0$ has been chosen so that
the ball $B_R:=\{ X \in \mathfrak t\mid \|X\|\leq R\}$ is contained in
the set
\begin{equation}
\label{eq:S}
S:=\{ X=iH \in \mathfrak t\mid \text{$|\alpha(H)|\leq \pi/2$ for all $\alpha\in \Delta$}\}\,.
\end{equation}
Then the map $\exp:\mathfrak t \to T$ is a diffeomorphism of $B_R$ onto its image,
say $D_R$. We shall refer to this condition by saying that $R$ is \textit{small}. We set
$C^\infty_R(T)^W:=\{f \in C^\infty(T)^W\mid \Supp(f) \subseteq D_R\}$.
Moreover, we denote by $PW_R(\mathfrak t^*)$ the space of functions on $P^+$ admitting
an extension $F:\mathfrak a_\C^*=\mathfrak t_\C^* \to \C$
so that $\lambda \mapsto \bc(m,\lambda)F(\lambda-\rho)$
is a $W$-invariant holomorphic function of exponential type $R$
(see \cite{BOP1}, Definition 3.2).

The following theorem is the main result of this note.

\begin{theorem}[Local Paley-Wiener theorem for the Jacobi transform]
\label{thm:PW}
Suppose that the root system is reduced, all multiplicities $m_\alpha$
are even and that $R>0$ is small.
Then the Jacobi transform is a bijection of $C^\infty_R(T)^W$ onto
$PW_R(\mathfrak t^*)$.
\end{theorem}

The remaining of this section is devoted to the proof of Theorem \ref{thm:PW}.

To show that the Jacobi transform maps into $PW_R(\mathfrak t^*)$ one can
follow the same lines as in Theorem 3.8 of \cite{BOP1}. We only remark here
that the exponential type $R$ for $\bc(m,\lambda)\widehat{f}(m,\lambda-\rho)$
depends on (\ref{eq:Jacobiexpl}) together with the fact that $d(m,\lambda)$
is a polynomial function (cf.\ also \cite{BOP1}, Lemma 3.3).

To prove the surjectivity, let $F \in PW_R(\mathfrak t^*)$ and set, according to
the inversion formula (\ref{eq:invJacobi}),
\begin{equation}
\label{eq:inverse}
f(t):=\sum_{\mu \in P^+} F(\mu) \dfrac{P(m,\mu,t)}{\|P(m,\mu)\|^2_m}\,.
\end{equation}

We need to show:
\begin{enumerate}
\item $f$ is smooth and $W$-invariant;
\item $\widehat{f}=F$;
\item $\Supp (f) \subseteq B_R$.
\end{enumerate}

For each $\alpha \in \Delta$, let $A_\alpha\in \mathfrak a$ be defined by
the condition that $\alpha(H)=\inner{A_\alpha}{H}$ for all $H \in \mathfrak a$.
Then $\Gamma:={\rm span}_\Z \{2\pi i A_\alpha/\inner{\alpha}{\alpha}\mid \alpha \in \Delta^+\}$
is a lattice in $\mathfrak t$ so that $A_\C=\mathfrak a_\C/\Gamma$.

\begin{lemma}\label{lemma:invEucl}
Let $R>0$ be small according to the definition given at the beginning of this
section. For $F \in PW_R(\mathfrak t^*)$ define $h:\mathfrak t \to \C$ by
$$h(H):=\sum_{\mu\in P} \bc(m,\mu)F(\mu-\rho) e^{\mu(H)}\,.$$
Then $h$ is a $\Gamma$-periodic function on $\mathfrak t$ and $\Supp (h) \subseteq
B_R + \Gamma$.
\end{lemma}
\begin{proof}
See \cite{BOP1}, Lemma 3.14.
\end{proof}

\begin{lemma}\label{lemma:PW}
Let $R>0$ be small and let $S$ be as in \nnn{eq:S}. Let $F \in PW_R(\mathfrak t^*)$
and define $f$ by \nnn{eq:inverse}. Then the following properties hold:
\begin{enumerate}
\item $f$ is smooth and $W$-invariant on $S$.
\item For all $t \in T$ we have, with $\widetilde{C}$ as in \nnn{eq:Ctilde},
\begin{equation}\label{eq:invwithD}
\delta(m,t)f(t)=\widetilde{C}
\, D\left(\sum_{\mu\in P} \bc(m,\mu) F(\mu-\rho)t^{\mu}\right).
\end{equation}
\item $f$ extends as a $W$-invariant smooth function on $T$ so that
$\widehat{f}=F$ and $\Supp (f) \subseteq D_R$.
\end{enumerate}
\end{lemma}

\begin{proof}
Observe first that ${\|P(m,\mu)\|^{-2}_m}$ is uniformly bounded in $\mu \in P^+$.
According to Lemma \ref{lemma:estJacobi}, for each multi-index $I$, each $t \in S$ and
each $N\in\N$, there is a constant
$C_{t,I,N}>0$ such that
$$\left| F(\mu) \, \dfrac{\partial_I P(m,\mu,t)}{\|P(m,\mu)\|_m^2}\right|
\leq C_{t,I,N} (1+\|\mu\|)^{-N+|I|+\sum_{\alpha\in\Delta^+} m_\alpha/2}
\,.$$
By choosing $N$ large enough, it follows that the series
$\sum_\mu F(\mu) \usfrac{\partial_I P(m,\mu)}{\|P(m,\mu)\|^2}$ converges uniformly. Thus $f$
is smooth. This proves (1).

The set $P^+$ is a fundamental domain for the action of $W$ on $P$, and $\rho \in P$
since $m$ is even. Because of (\ref{eq:Dzero}), we obtain as in \cite{BOP1}, Lemma 3.13:
$$
D\left(\sum_{\mu\in P^+,\,w\in W} \bc(m,\mu+\rho)F(\mu) t^{w(\mu+\rho)}\right)
=D\left(\sum_{\mu\in P} \bc(m,\mu)F(\mu-\rho) t^{\mu}\right)\,.
$$
It follows by Lemma \ref{lemma:invEucl} that
\begin{align*}
\delta(m,t) f(t)&=\sum_{\mu\in P^+} F(\mu) \frac{\delta(m,t) P(m,\mu,t)}{\|P(m,\mu)\|^2_m} \\
&=\widetilde{C}\; \bc(m,\mu+\rho)
\sum_{\mu\in P^+} F(\mu) D\left(\sum_{w \in W} t^{w(\mu+\rho)}\right) \\
&=\widetilde{C} \;
 D\left(\sum_{\mu\in P} \bc(m,\mu)F(\mu-\rho) t^{\mu}\right).
\end{align*}
By Lemma \ref{lemma:invEucl}, the right-hand side of (\ref{eq:invwithD})
is supported in $D_R$. Thus $f$ extends as smooth $W$-invariant function
on $T$ with $\Supp (f) \subseteq D_R$. Finally, as $f \in L^2(T)^W$, we have
$$\sum_{\mu\in P^+} F(\mu) \frac{P(m,\mu,t)}{\|P(m,\mu)\|^2_m}
=f(t)=\sum_{\mu\in P^+} \widehat{f}(m,\mu) \frac{P(m,\mu,t)}{\|P(m,\mu)\|^2_m}\,$$
and hence $F(\mu)=\widehat{f}(m,\mu)$ for all $\mu\in P^+$.
\end{proof}

As a corollary of the proof of the Paley-Wiener theorem, we obtain two integral formulas for functions in $C^\infty_R(T)^W$ for $R>0$ small and
even multiplicities.

\begin{lemma}\label{lem:intformula}
Suppose that $R>0$ is small and that $f\in C^\infty_R(T)^W$. Then
\begin{equation}
\label{eq:intformula}
\delta(m,t)f(t)=\widetilde{C} \; D\left( \int_{i\mathfrak t^*} \bc(m,\lambda)
\widehat{f}(m,\lambda-\rho) t^\lambda \; d\lambda\right)\,,
\end{equation}
where $D$ is the differential operator of (\ref{eq:formulaF}) and $\widetilde{C}$ is
as in \nnn{eq:Ctilde}. Moreover,
$$
f(t)=\widetilde{C} \int_{i\mathfrak t^*} \widehat{f}(m,\lambda-\rho) F(m,\lambda,t)
e(m,\lambda) \; d\lambda\,,$$
where
$$e(m,\lambda):=\dfrac{1}{C}\; \prod_{\alpha\in\Delta^+} \prod_{k=0}^{m_\alpha/2-1}
(\lambda_\alpha-k)\,.$$
\end{lemma}
\begin{proof}
As Corollary 3.16 of \cite{BOP1}, using (\ref{eq:invwithD}) and (\ref{eq:formulaF}).
\end{proof}

\section{The local Huygens' principle for the modified Laplace operator
on $T$ in the even multiplicity case}
\noindent
Let $m\geq 0$ be fixed. The Laplace operator on the
compact torus $T$ is the $W$-invariant differential operator
$$L(m)=L_T - \sum_{\alpha\in \Delta^+} m_\alpha
\dfrac{1+e^{-2\alpha}}{1-e^{-2\alpha}} \, \partial_{i\alpha}\, $$
where $L_{\mathfrak t}$ is the Laplace operator on
the abelian Lie group $T$ and
$\partial_{i\alpha}$ denotes the derivative in the direction of
the vector $iA_\alpha$. Note that the second
term vanishes if $\alpha(H)\in i\frac{\pi}{2}+i\pi\Z$ for all
$\alpha$ and has singularities if there
is a root $\alpha$ such that $\alpha (H)\in i\pi\Z$.

The Jacobi polynomials are eigenfunctions of
$L(m)$, with $L(m)P(m,\mu)=-\inner{\mu+2\rho}{\mu}P(m,\mu)$ for all
$\mu \in P^+$.
The \textit{modified wave equation} on $T$ is the partial differential equation
\begin{equation}\label{eq:waveeq}
\big(L(m)-\|\rho\|^2\big)u=u_{\tau\tau}
\end{equation}
where $u=u(t,\tau)$ is a function of $(t,\tau)\in T \times I$ and
$I \subseteq \R$ is an interval containing $0$.
As far as we know, there is no statement in the literature about
finite propagation speed for the solutions of this wave equation.
This will follow from the following lemma; cf.\ Lemma 4.3 in \cite{BOP1}:

\begin{lemma}
Let $T$ be the torus associated with a triple
$(\mathfrak a, \Delta, m)$ with reduced root system $\Delta$ and
even multiplicity function $m$.
Let $R>0$ be small according to the previous section. Let
$0<\e<R$, and let $f \in C^\infty_\e (T)^W$. Assume that
$u(t,\tau)$ is a distributional solution to the Cauchy problem
\begin{equation}\label{eq:CauchyT1}
\left\{\begin{array}{rl}
(L(m)-\|\rho\|^2)u\!\!\!\!\!&=u_{\tau\tau}\,, \\
u(t,0)\!\!\!\!\!&=0\,,  \\
u_\tau(t,0)\!\!\!\!\!&=f(t)\,.
\end{array}\right.
\end{equation}
We suppose that $u$ is smooth and $W$-invariant in the variable
$t\in T$. Then for $t=\exp (X) \in T$ and $\tau\in [0,R-\e]$ we have
\begin{equation}\label{eq:udiff}
\delta (m,t) u(t,\tau)= \widetilde{C}\,  D \int_{i\mathfrak t^*}\bc (m,\lambda)  \widehat{f}(m,\lambda-\rho)
\frac{\sin ( \|\lambda \| \tau )}{\|\lambda \|} t^\lambda \; d\lambda\, ,
\end{equation}
where $D$ is the differential operator of (\ref{eq:formulaF}) and $\widetilde{C}$ is
as in \nnn{eq:Ctilde}.

In particular the following hold for $\tau\in [0,R-\e]$:
\begin{enumerate}
\item $u$ is a smooth function of $\tau$;
\item  \emph{{\rm(}Finite propagation speed{\rm)}} $u(\exp X,\tau )=0$ if
$\|X\|\ge |\tau |+\e$;
\item \emph{{\rm(}Local strong Huygens' principle{\rm)}} $u(\exp X,\tau )=0$ if $\mathfrak{t}$ is
odd dimensional and $|\tau |\ge \|X\|+\epsilon$.
\end{enumerate}
\end{lemma}
\begin{proof}
Since the domain of integration $T$ for the Jacobi transform is compact, we can interchange differentiation and integration.
Applying the Jacobi transform in the $t$-variable to both sides
of (\ref{eq:CauchyT1}) implies therefore the
following initial value problem in the
$\tau$-variable for $\widehat{u}(m,\mu,\tau)$ :
\begin{equation}\label{eq:CauchyT2}
\left\{\begin{array}{rl}
\widehat{u}_{\tau\tau}(m,\mu,\tau)\!\!\!\!\!&=-\|\mu+\rho\|^2\widehat{u}(m,\mu,\tau) \\
\widehat{u}(m,\mu,0)\!\!\!\!\!&=0 \\
\widehat{u}_{\tau}(m,\mu,0)\!\!\!\!\!&=\widehat{f}(m,\mu)\, .
\end{array}\right.
\end{equation}
Thus
$$\widehat{u}(m,\mu,\tau)=\widehat{f}(m,\mu)\,
\dfrac{\sin (\|\mu +\rho\|\tau )}{\|\mu +\rho\|}\, .$$
The local Paley-Wiener theorem now implies that $\mu \mapsto
\widehat{f}(m,\mu)$ belongs to $PW_\e(\mathfrak t^*)$.
The function $\mu \mapsto \usfrac{\sin (\|\mu \|\tau) }{\|\mu \|}$
has also a $W$-invariant holomorphic extension of exponential type
$\tau$. Thus the function $\mu \mapsto \bc(m,\mu) \widehat{f}(m,\mu - \rho )\, \usfrac{\sin (\|\mu \|\tau) }{\|\mu \|}$ admits a $W$-invariant holomorphic extension
of exponential type $\e+\tau$, that is $\widehat{u}(m,\mu,\tau)=\widehat{f}(m,\mu)
\usfrac{\sin (\|\mu +\rho\|\tau ) }{\|\mu +\rho\|}$ belongs to
$PW_{\e+\tau}(\mathfrak t^*)$. Since $\e+\tau$ is small for $\tau \in [0,R-\e]$, another application of the local Paley-Wiener theorem yields $u(\cdot,\tau)\in C^\infty_{\e+\tau}(T)^W$. This proves the finite propagation speed and allows us to apply Lemma \ref{lem:intformula} to $u(\cdot,\tau)$. Thus (\ref{eq:udiff}) holds.

Finally, recall from (4.13) in \cite{BOP1} that the function
$$v(X,\tau):=\int_{i\mathfrak{t}^*}\bc (m,\lambda)
\widehat{f}(m,\lambda-\rho )
\frac{\sin (\|\lambda \|\tau )}{\|\lambda \|} e^{\lambda(X)}\, d\lambda$$
is the solution of the Cauchy problem for the wave equation on
$\mathfrak t \cong \R^n$:
$$
\left\{\begin{array}{rl}
L_{\mathfrak{t}}v(X,\tau)\!\!\!\!\!&=v_{\tau\tau} (X,\tau )\\
v(X,0)\!\!\!\!\!&=0\\
v_\tau(X,0)\!\!\!\!\!&=g(X)\,,
\end{array}\right.
$$
where $g\in C_\e^\infty(T)^W$ is the inverse Fourier transform of $\lambda \mapsto \bc (m,\lambda)\widehat{f}(m,\lambda -\rho)$.
By the standard results for the
wave equation on Euclidean spaces it follows that, in case $\mathfrak{k}$ is odd dimensional, $v(X,\tau)=0$
for $|\tau |\ge \|X\|+\epsilon$ (the strong Huygens principle).
Let $U$ be the periodization of $v$ in the $X$-variable, and note
that the translates of the support of $g$ are all
disjoint. Then (\ref{eq:udiff}) implies that
that $u(\exp X,\tau )=\widetilde{C} D U(X,\tau )$ safisfies (3).
\end{proof}

Different proofs of the local strong Huygens principle can be obtained
from the results in the previous section together with same arguments employed in \cite{BOP1}. More precisely, one can prove the following theorem.

\begin{theorem}\label{thm:HP}
Let $T$ be the torus associated with a triple $(\mathfrak a, \Delta, m)$ with reduced root system $\Delta$ and
even multiplicity function $m$.
Let $R>0$ be small according to the previous section. Let
$0<\e<R$, and let $f \in C^\infty_\e(T)^W$.
Assume that $T$ is odd dimensional.

Suppose moreover that $u(t,\tau)$ is a smooth solution of the
Cauchy problem
\begin{equation}\label{eq:CauchyT}
\left\{\begin{array}{rl}
(L(m)-\|\rho\|^2)u\!\!\!\!\!&=u_{\tau\tau}\,, \\
u(t,0)\!\!\!\!\!&=0\,,  \\
u_\tau(t,0)\!\!\!\!\!&=f(t)\,.
\end{array}\right.
\end{equation}
Then the following properties are satisfied:
\begin{enumerate}
\renewcommand{\theenumi}{\alph{enumi}}
\renewcommand{\labelenumi}{(\theenumi)}
\item \emph{(Local exponential Huygens' principle)}
There is a constant $C>0$ so that
for all $(t,\tau) \in T\times [0,R-\e]$
and all $\gamma \in [0,\infty)$ we have
\begin{equation} \label{eq:expHP}
|\delta(m,t)u(t,\tau)|\leq C e^{-\gamma(\tau-\|H\|-\e)}\,
\end{equation}
with $t=\exp H$.
\item \emph{(Local strong Huygens' principle)}
$$\Supp(u) \cap ( T \times [0,R-\e])=
\Supp(u) \cap ( D_R \times [0,R-\e]) \subseteq S_\e\,,$$
where $S_\e$ denotes the $\e$-shell
\begin{equation} \label{eq:eshell}
S_\e:=\{(t=\exp H,\tau) \in T \times [0,\infty)\mid \tau-\e\leq \|H\|\leq \tau+\e\}.
\end{equation}
\item
Suppose $\dim(T) \geq 3$. Let $D$ be the differential operator of
\nnn{eq:formulaF}.
Then for all $t=\exp H \in T$ and $\tau \in [0,R-\e]$
the smooth solution $u(t,\tau)$ to \straightparen{\ref{eq:CauchyT}} is given by the formula
\begin{equation}
\label{eq:uodd}
\delta(m,t)u(t,\tau)=\frac{\Omega_n/2}{[(n-3)/2]! \Omega_{n-1}}\;
D\left( \frac{\partial}{\partial(\tau^2)}\right)^{(n-3)/2}
  \big( \tau^{n-2} (M^\tau g)(H)\big)\,.
\end{equation}
Here $g\in C^\infty_\e(\mathfrak t)^W$ is the inverse Euclidean Fourier
transform
of the function $\lambda \mapsto \bc(m,\lambda)\widehat{f}(m,\lambda-\rho)$, and
\begin{equation}\label{eq:M}
(M^rg)(H):=\frac{1}{\Omega_{n-1}(r)} \int_{S_r(H)} g(s)\,d\sigma(s)\,
\end{equation}
is the mean value of a function $g:\mathfrak t\to \C$ on the Euclidean sphere
  $S_r(H):=\{H\in \mathfrak t\mid
\|H\|=r\}$ in $\mathfrak t \cong \R^n$ with respect to the
${\rm O}(n)$-invariant surface measure $d\sigma$. Moreover,
$\Omega_{n-1}(r)$ denotes the surface area of $S_r(H)$, and
$\Omega_{n-1}:=\Omega_{n-1}(1)$.
\end{enumerate}
\end{theorem}

\end{document}